\documentclass[12pt]{article}
\usepackage[]{geometry}

\usepackage{latexsym}
\usepackage{amssymb}
\usepackage{amsmath}
\usepackage{amsthm}
\usepackage{amscd}

\usepackage[pdftex]{graphicx}

\newtheorem{thm}{Theorem}[section]
\newtheorem{cor}[thm]{Corollary} 
\newtheorem{lem}[thm]{Lemma} 
\newtheorem{prop}[thm]{Proposition} 
\newtheorem{rem}[thm]{Remark} 
\newtheorem{exam}[thm]{Example} 
 
\newtheorem{defin}[thm]{Definition}

\newcommand{\demo}{\par\noindent{\it Proof. \/}\ }
\newcommand{\enD}{\hfill $\Box$\vspace{3truemm} \par}
\newcommand{\R}{\mathbb{R}}

\newcommand{\sign}{\operatorname{sign}}

\begin{document}
\title{A view-parametric extension of the d'Ocagne formula for a surface in $\R^3$
}

\author{Ken Anjyo and Yutaro Kabata}

\date{\today}

\maketitle
\begin{abstract}
\noindent In this paper, we consider the orthogonal projection of a surface in $\R^3$ for a given view direction. We then introduce and investigate several invariants of the families of the plane curves that locally configure the projection image of the surface. Using the invariants, we also show an extension of the d'Ocagne formula that associates a local behavior of the projection image of a surface with Gaussian curvature of the surface.
\end{abstract}

\renewcommand{\thefootnote}{\fnsymbol{footnote}}

\section{Introduction}
Jan Koenderink showed in his famous book \cite{Koenderink} that classical differential geometry is a strong tool to study vision science.
Among his various results in the area, what is called Koenderink's formula (see \cite{Koenderink1984, Koenderink}) has attracted many mathematicians \cite{Araujo, Cipolla-Giblin, Damon-Giblin-Haslinger, FHS, HKS, IRFT}.
This formula gives the relationship between the Gaussian curvature of a surface and the curvature of the apparent contour.
For instance, while it is a naive intuition that we know whether a surface is locally elliptic or hyperbolic from its apparent contour, the formula mathematically justifies this intuition.

A precise statement of the Koenderink formula is given as follows:
Let $M \subset \R^3$ be a smooth surface and $K$ the Gaussian curvature of $M$.
For a point $p\in M$ and a tangent vector $v\in T_p M$,
let $\kappa_n(p)$ denote the normal curvature of $M$ at $p$ along $v$.
When $\kappa_n(p)\not=0$, the apparent contour (the singular value set of the orthogonal projection of $M$ along $v$) is non-singular at $p$. Denoting the curvature of the contour as a plane curve by $\kappa_c(p)$, we have
the Koenderink formula:
\begin{equation}\label{KoenderinkFormula}
    \kappa_c(p) \kappa_n(p)=K(p).
\end{equation}
Recently, it is found that the same formula was already discovered by d'Ocagne in 1895 (see also \cite[page 22]{Araujo}). We therefore refer to the equation (\ref{KoenderinkFormula}) as {\it the d'Ocagne formula} throughout this paper.
{Note that the above expression of the formula in terms of curvatures in surface theory was formulated by Koenderink \cite{Koenderink}.}

\begin{figure}
\centering
\includegraphics[clip,width=8.0cm]{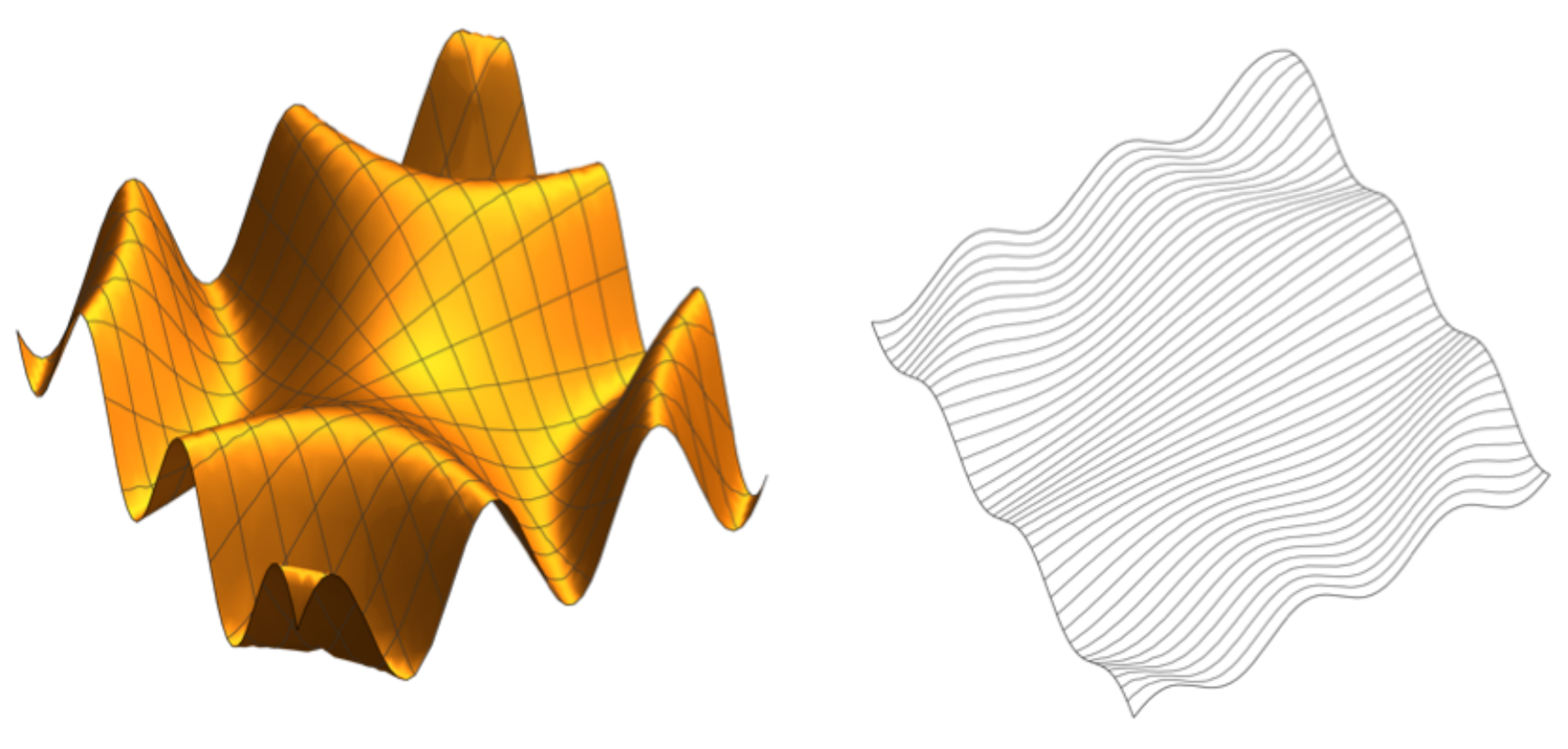}
\caption{A surface given by $z=\sin xy$. The right side is just drawn by the curves parametrized by the parameter $y$, that is, each space curve is expressed as $(x_0,y,\sin x_0 y)$ for a fixed $x_0$.}
\label{curvessurffig}
\end{figure}

While Koenderink in \cite{Koenderink} is interested mainly in the boundary of the projection image (apparent contour) of a surface, we are interested in the interior of its projection image.
Suppose that a given surface consists of a certain family of space curves, then its projection image gives a one-parameter family of curves as illustrated in Figure \ref{curvessurffig}.
This situation would make us expect 
that there is a geometric relationship between the shape of the surface and the pattern of the projection image.
The present paper aims to justify this expectation by deriving formulae similar to the d'Ocagne formula. To achieve this goal, assuming that all maps and surfaces considered in this paper are smooth ($C^\infty$), we investigate the above one-parameter family of curves in the projective image of a surface and then recovering the Gaussian curvature of the surface.
Differential geometric approaches for families of curves as in \cite{Kabata-Takahashi, Takahashi2} would then be quite useful.

\begin{figure}
\centering
\includegraphics[clip,width=12.0cm]{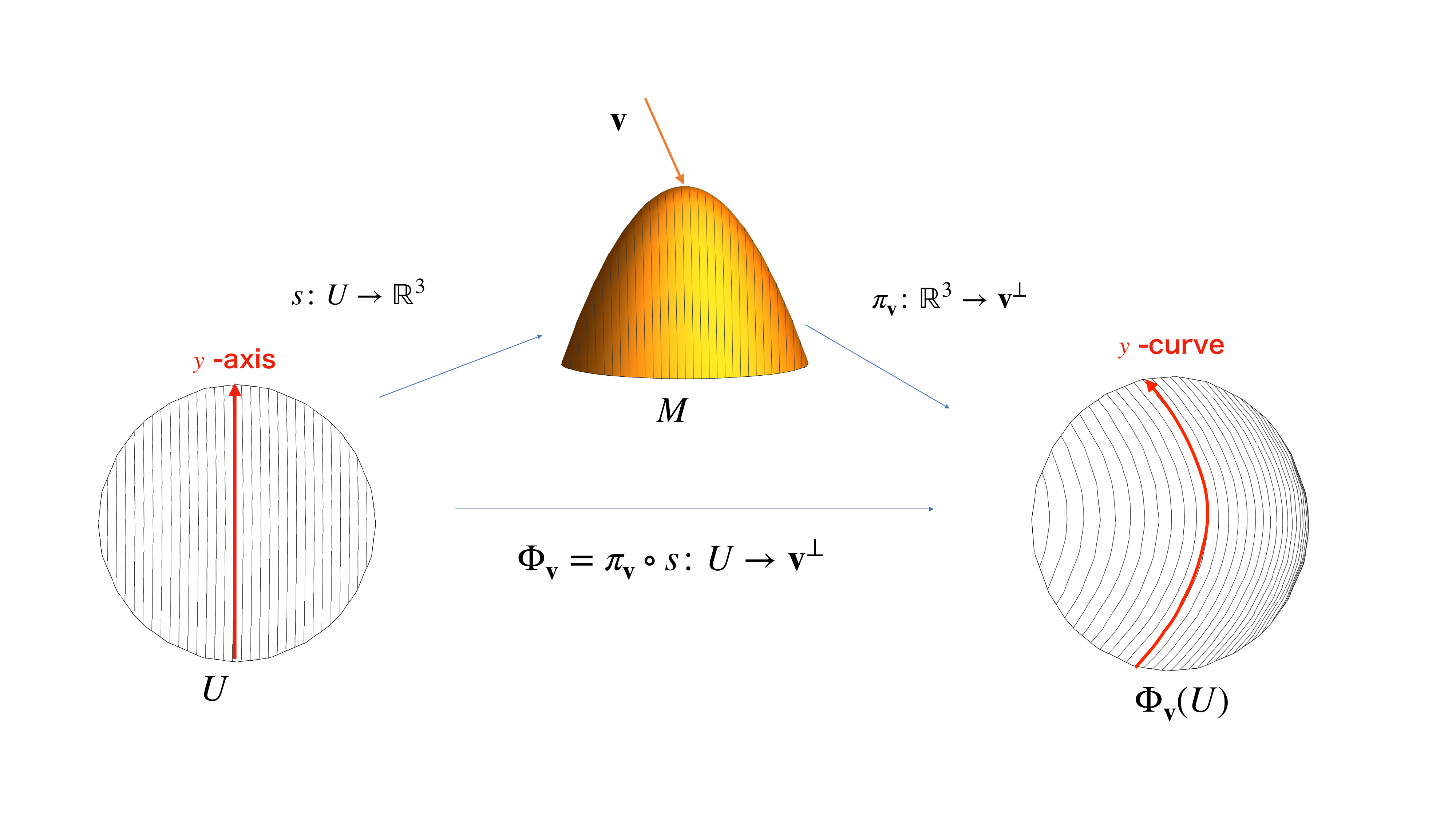}
\caption{The orthogonal projection of a surface $M$ along the view direction ${\bf v}$ and the projection image giving a one-parameter family of curves.
}
\label{fig:ellipproj}
\end{figure}

The rest of this paper is organized as follows. As a preliminary step, we define several differential geometric notions for {one-parameter families of regular curves in plane} in \S 2, which is important in stating our main results. In particular, we note that these differential geometric notions
are invariant under translations and rotations on the plane (see Remark \ref{rem:invariant}).

In \S 3, we consider orthogonal projections of a surface in $\R^3$ that consists of certain families of space curves.
Here we suppose that the surface is expressed as the graph of a function $g(x,y)$: $s(x,y)=(x,y,g(x,y))^t$, where $(x, y)^t$ belongs to a certain domain $U$ in $\R^2$.
For a fixed $x$ (resp. $y$), $s(x,y)$ gives a space curve parametrized by $y$ (resp. $x$), thus our surface is regarded as a one-parameter family of space curves $\{s(x_0,y) \}_{x_0}$ (resp. $\{s(x,y_0) \}_{y_0}$).
The orthogonal projection of such a surface is then considered for a given unit vector $\mathbf{v}$, which is referred to as the \textit{view direction} in this paper. As shown in Figure 2, it means that we investigate the projection $\pi_{\mathbf{v}}$ of the surface $s$ onto a plane $\mathbf{v}^{\perp}$ which is perpendicular to the view direction $\mathbf{v}$. More specifically we examine the
one-parameter families of plane curves on $\mathbf{\Phi}_{\mathbf{v}}(U)$, where $\mathbf{\Phi}_{\mathbf{v}}:= \pi_{\mathbf{v}}\circ s$.  
Thus we derive several formulae which describe relationship between the invariants of such one-parameter families of plane curves and geometric information of the surface. 

In \S 4, we investigate the signs of the formulae obtained in \S 3. This leads us to consider an extension of the d'Ocagne formula. Actually, in a special setting in \S 5, the formulae containing the Gaussian curvature are given. We note all the formulae obtained in this paper include the two parameters that prescribe a view direction. In this sense, we say those formulae are expressed along with the view direction parameters.


\bigskip
\noindent
{\bf Acknowledgement}. 
The second author is partially supported by JSPS KAKENHI Grant Number JP 20K14312. 
The authors want to thank Professor Farid Tari for his helpful comment.


\section{One-parameter family of smooth curves and its invariant}
We consider local differential geometry of families of plane curves,
which is regarded as a smooth mapping from plane to plane.
We assume that $\R^2$ is the Euclidean plane equipped with the scalar product $<\cdot,\cdot>$, and denote the Euclidean distance by $||\cdot||$.
Let $U\subset \R^2$ be a simply connected open subset.
\begin{defin}\label{def:oneparacyrve}
{\rm
We say that a smooth map $f\colon U \to \R^2, (u,v)\mapsto f(u,v)$ is {\it a one-parameter family of regular curves with respect to $u$ (resp. $v$)}
when $f_u(u,v)\not=0$ (resp. $f_v(u,v)\not=0$) for all $(u,v)\in U$.
}
\end{defin}

Let $f\colon U \to \R^2, (u,v)\mapsto f(u,v)$ be a one-parameter family of regular curves with respect to $u$ (resp. $v$).
For a fixed parameter $v$ (resp. $u$),
$f(\cdot,v)$ (resp. $f(u,\cdot)$) is a regular plane curve,
which is called {\it a $u$ (resp. $v$)-curve}.
Such $u$ or $ v$-curves have the curvatures in the usual sense.
\begin{defin}\label{def:curvature}
{\rm
We define {\it the curvature of a one-parameter family of regular curves $f$ with respect to $u$ (resp. $v$)}
as the function
$$
\kappa[f,u](u,v):=\frac{\det (f_u\; f_{uu})}{\|f_u\|^{3}}(u,v)
\;\; \left(\mbox{\it resp.}\;\; \kappa[f,v](u,v):=\frac{\det (f_v\; f_{vv})}{\|f_v\|^{3}}(u,v)\right).
$$
Here $\det (\mathbf{a}\; \mathbf{b})$ for $\mathbf{a}, \mathbf{b}\in\R^2$ means the determinant of the matrix whose columns are $\mathbf{a}$ and $\mathbf{b}$.
}
\end{defin}
\begin{exam}
    Figure \ref{exuv2}  shows the image of the mapping $f_\pm(u,v)=(u\pm v^2,v)$ which are regarded as one-parameter families of regular curves with respect to $v$ with $\kappa[f_\pm,v]\lessgtr 0$.
\begin{figure}
\centering
  \includegraphics[width=6.0cm]{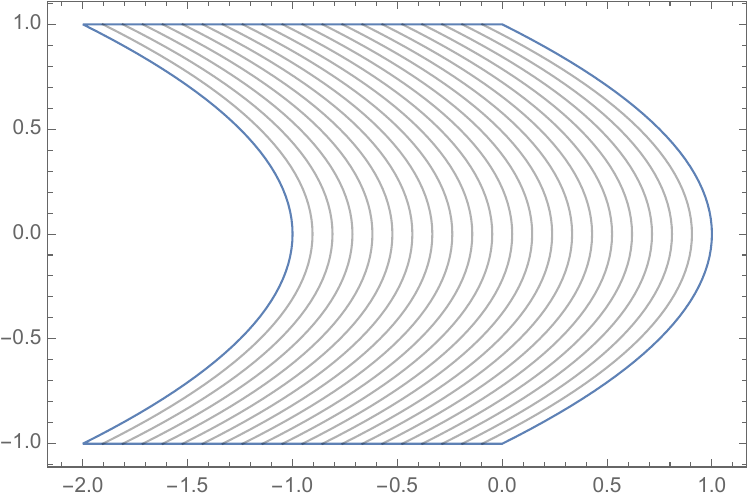}
  \hspace{2cm}
  \includegraphics[width=6.0cm]{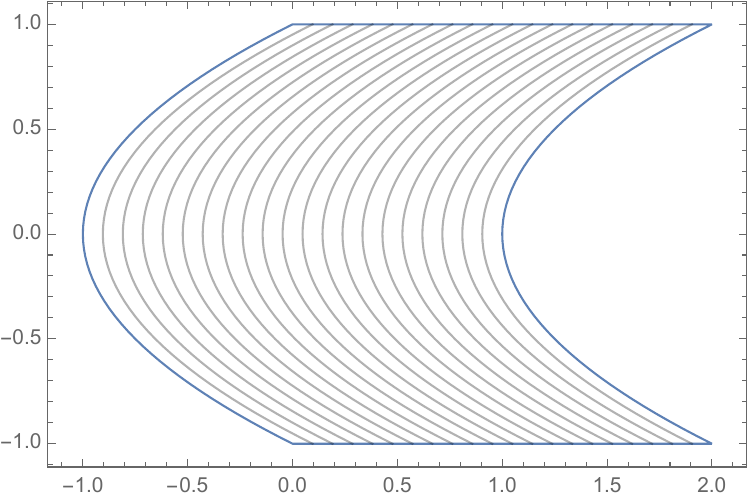}
\caption{The left (resp. right) represents the image of $f_-$ (resp. $f_+$) for $f_\pm\colon (-1,1)\times (-1,1)\to \R^2$, $f_\pm(u,v)=(u\pm v^2,v)$.} 
 \label{exuv2} 
\end{figure}

\end{exam}


Next, we define the squared velocity of a smooth mapping.
\begin{defin}\label{def:sqvel}
{\rm
Let $f\colon U \to \R^2, (u,v)\mapsto f(u,v)$ be a smooth mapping.
We define {\it the squared velocity $SV[f,u](u,v)$ (resp. $SV[f,v](u,v)$) of $f$ at $(u,v)\in U$ with respect to $u$ (resp. $v$)} as 
$$
SV[f,u](u,v):=\|f_u(u,v)\|^2 \;\; 
(\mbox{\it resp.}\;\;  SV[f,v](u,v):=\|f_v(u,v)\|^2).
$$
}
\end{defin}

\begin{figure}[htbp]
\centering
  \includegraphics[width=10.0cm]{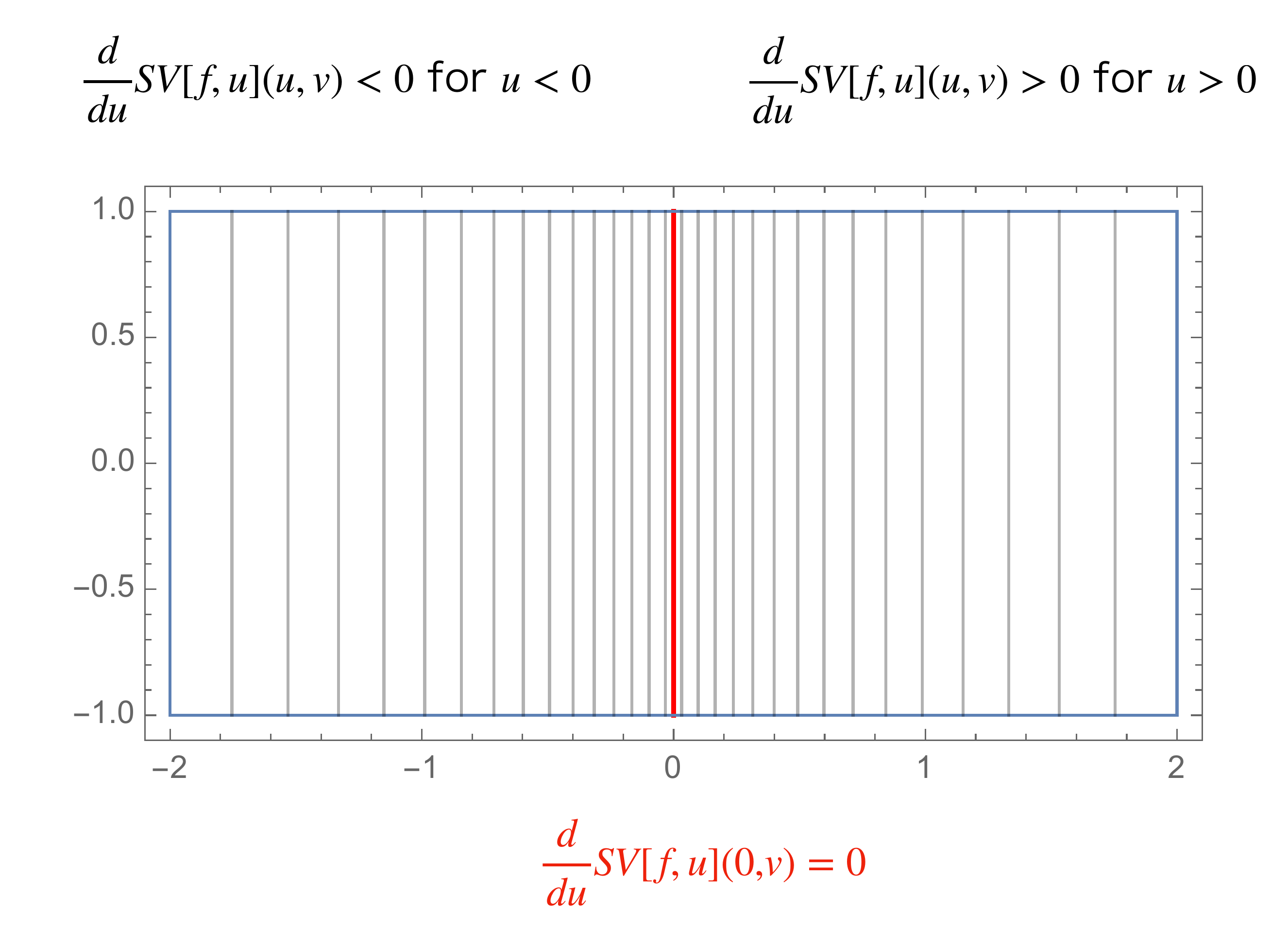}
\caption{The image of the mapping $f\colon (-1,1)\times (-1,1)\to \R^2$, $f(u,v)=(u+ u^3,v)$. The red line is the image of the $v$-axis in the $uv$-plane by $f$.} 
 \label{exu+u^3} 
\end{figure}
\begin{exam}
    Figure \ref{exu+u^3} shows the image of the mapping $f(u,v)=(u+u^3,v)$ which is regarded as a one-parameter family of curves (straight lines) parametrized by $v$. Here we have $SV[f,u]=(1+3u^2)^2$ and $\frac{d}{du} SV[f,u]=12u(1+3u^2)$, thus $\frac{d}{du} SV[f,u] <0$ for $u<0$; $\frac{d}{du} SV[f,u] >0$ for $u>0$; and $\frac{d}{du} SV[f,u] =0$ for $u=0$.
\end{exam}

\begin{rem}{\rm
Let $f, \tilde{f} \colon U \to \R^2$ be smooth mappings such that $\tilde{f}(u,v)=A(f(u,v))+a$ for a rotation $A\in SO(2)$ and a constant vector $a\in\R^2$. Then $SV[f,u]=SV[\tilde{f},u]$ and $SV[f,v]=SV[\tilde{f},v]$ hold. In addition, assuming that $f$ and  $\tilde{f}$ are one-parameter families of regular curves with respect to $u$ (resp. $v$), then $\kappa[f,u]=\kappa[\tilde{f},u]$ (resp.  $\kappa[f,v]=\kappa[\tilde{f},v]$) holds.
In this sense, the differential geometric notions introduced in this section are invariant under the action of rotations and translations on the plane. In the following sections, the partial differentials of the squared velocities: $\frac{d}{du} SV[f,u]$ and $\frac{d}{dv} SV[f,v]$ also play important roles and are invariant in the above sense.
Note that a general theory for invariants of families of plane curves is considered in \cite{Kabata-Takahashi}. 
In particular, the notion of one-parameter families of Legendre curves is introduced, and their complete invariants under the action of rotations and translations on the plane are given with the uniqueness and existence theorem.
}
\label{rem:invariant}
    \end{rem}



\section{Formulae with respect to projections of a surface}
\begin{figure}[htbp]
\centering
  \includegraphics[width=14.0cm, clip]{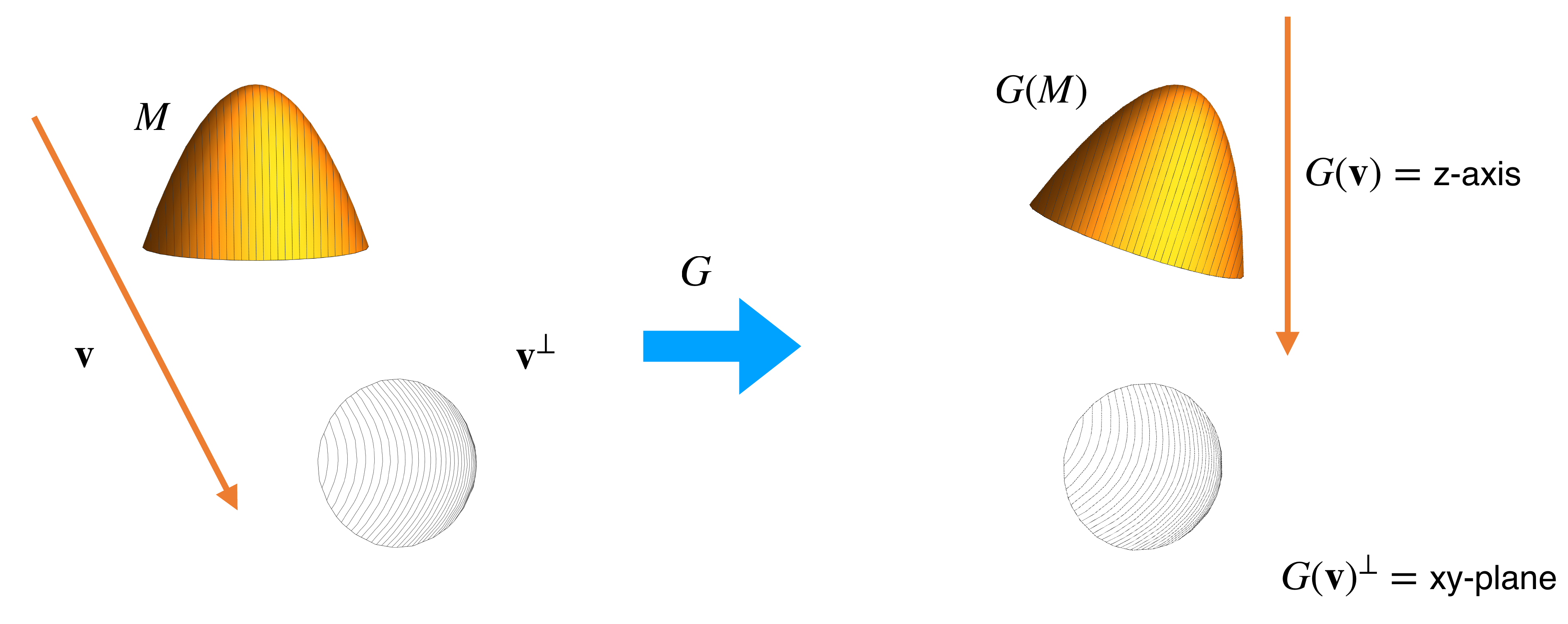}
\caption{The orthogonal projection of a surface $M$ along a view direction ${\bf v}$ and its rotation.} 
 \label{fig:projandrot} 
\end{figure}

We consider local geometry of a smooth surface $M\subset \R^3$ through projection.
Let $\pi_{\bf v}\colon\R^3 \to {\bf v}^\bot$ be the orthogonal projection with the view direction ${\bf v}\in S^2$.
We consider the family of regular curves generated by the orthogonal projection $\pi_{\bf v}|S$
as in Figure \ref{fig:ellipproj}--\ref{fig:paraboproj}.
For convenience,
we rotate our surface $M$ so that ${\bf v}$ is moved into $(0,0,-1)^t$,
and then our orthogonal projection is regarded as the canonical projection from $\R^3$ to the $xy$-plane. See Figure \ref{fig:projandrot}.

Without loss of generality, we may suppose that $M$ contains the origin $0\in\R^3$,
and $M$ around $0$ is locally parameterized by $s\colon U\to \R^3, \;s(x,y)=(x,y,g(x,y))^t$ for $(x,y)^t\in U$
with a simply connected open subset $U\subset\R^2$ containing the origin $0\in\R^2$
and a smooth function $g\colon U\to \R$ so that $g(0)=0$.
Take the view direction ${\bf v}=(\sin \theta \cos \phi,\sin \theta \sin \phi, \cos \theta)^t\in S^2$
with $\theta, \phi \in \R$.
Set the rotation matrices as
$$
R_y(a)=\left(
\begin{array}{ccc}
\cos a & 0& \sin a \\
0&1&0\\
-\sin a & 0& \cos a \\
\end{array}
\right),\;\;\;\;
R_z(a)=\left(
\begin{array}{ccc}
\cos a& -\sin a & 0 \\
\sin a & \cos a& 0 \\
0&0&1\\
\end{array}
\right)
$$
with $a\in\R$.
Putting
$$
G:=R_y(\pi-\theta)R_z(-\phi)=\left(
\begin{array}{ccc}
-\cos\theta\cos\phi & -\cos\theta\sin\phi& \sin\theta \\
-\sin\phi&\cos\phi&0\\
-\sin\theta\cos\phi & -\sin\theta\sin\phi& -\cos\theta \\
\end{array}
\right),
$$
we have $G({\bf v})=(0,0,-1)^t$, and $G(M)=G(s(U))$ is parametrized as
$$
G(s(x,y))=
\left(
\begin{array}{c}
g(x,y) \sin \theta -x \cos \theta  \cos \phi -y \cos \theta  \sin \phi \\
y \cos \phi -x \sin \phi\\
-g(x,y) \cos \theta -x \sin \theta  \cos \phi -y \sin \theta  \sin \phi\\
\end{array}
\right).
$$
Thus the orthogonal projection of $G(M)=G(s(U))$ along $G({\bf v})$ is expressed as
$$
\Phi[{\bf v}](x,y)=
\left(
\begin{array}{c}
g(x,y) \sin \theta -x \cos \theta  \cos \phi -y \cos \theta  \sin \phi \\
y \cos \phi -x \sin \phi\end{array}
\right).
$$
Write
\begin{eqnarray}
    P_{\theta,\phi}(x,y):=-\cos\theta\cos\phi +g_x(x,y)\sin\theta,\\
    Q_{\theta,\phi}(x,y):=-\cos\theta\sin\phi +g_y(x,y)\sin\theta.
\end{eqnarray}
Immediately, we have the following lemma.

\begin{lem}
\begin{enumerate}
\item
$\Phi[{\bf v}]\colon U\to \R^2$ is  a one-parameter family of regular curves with respect to $x$  if and only if
either $\sin \phi\neq0$ or $P_{\theta,\phi}(x,y)\neq0$.
 \item
$\Phi[{\bf v}]\colon U\to \R^2$ is  a one-parameter family of regular curves with respect to $y$ if and only if
either $\cos \phi\neq0$ or $Q_{\theta,\phi}(x,y)\neq0$. 
\end{enumerate}
\label{regcondlem}
\end{lem}

With the condition in Lemma \ref{regcondlem},
$\Phi[{\bf v}]$ is a one-parameter family of regular curves with respect to $x$ or $y$.
Let $K$ denote the Gaussian curvature of the surface $M$. 
We have Theorems \ref{thm:y-zeta}, \ref{thm:x-zeta}. 

\begin{thm}\label{thm:y-zeta}
If $\Phi[{\bf v}]\colon U\to \R^2$ is  a one-parameter family of regular curves with respect to $y$,
then we have
$$
\kappa[\Phi[{\bf v}],y]
\cdot
\frac{d}{dx}\left( SV[\Phi[{\bf v}],x]\right)
=
-2g_{xx}g_{yy}\frac{P_{\theta,\phi}\cos\phi \sin^2\theta}{(\cos^2 \phi+Q_{\theta,\phi}^2)^{\frac32}}.
$$

\end{thm}

\begin{thm}\label{thm:x-zeta}
If $\Phi[{\bf v}]\colon U\to \R^2$ is  a one-parameter family of regular curves with respect to $x$,
then we have
$$
\kappa[\Phi[{\bf v}],x]
\cdot
\frac{d}{dy}\left( SV[\Phi[{\bf v}],y]\right)
=
2 g_{xx}g_{yy}\frac{Q_{\theta,\phi}\sin\phi \sin^2\theta}{(\sin^2 \phi+P_{\theta,\phi}^2)^{\frac32}}.
$$
\end{thm}

\noindent
{\it Proofs for Theorems \ref{thm:y-zeta}, \ref{thm:x-zeta}.}\;
 Combining statements in the following Propositions \ref{prop:curveturep}, \ref{prop:sv},
 we see the formulae.
\enD

\begin{prop}\label{prop:curveturep}
\begin{enumerate}
\item
If $\Phi[{\bf v}]\colon U\to \R^2$ is a one-parameter family of regular curves with respect to $y$,
then the curvature of the one-parameter family of regular curves $\Phi[{\bf v}]$ with respect to $y$ is
$$
\kappa[\Phi[{\bf v}],y]=
-\frac{g_{yy}\cos\phi \sin\theta}{(\cos^2 \phi+(\cos\theta\sin\phi-g_y\sin \theta)^2)^{3/2}}.
$$

\item
If $\Phi[{\bf v}]\colon U\to \R^2$ is  a one-parameter family of regular curves with respect to $x$,
then the curvature of the one-parameter family of regular curves $\Phi[{\bf v}]$ with respect to $x$ is
$$
\kappa[\Phi[{\bf v}],x]=
\frac{g_{xx}\sin\phi \sin\theta}{(\sin^2 \phi+(\cos\theta\cos\phi-g_x\sin \theta)^2)^{3/2}}.
$$
\end{enumerate}
\end{prop}
\demo
Direct calculations based on Definition \ref{def:curvature} show the statements.
\enD

\begin{prop}\label{prop:sv}
We have the following equations with respect to the differentials of the squared velocity of $\Phi[{\bf v}]$.
\begin{enumerate}
\item
$$
\frac{d}{dx}\left( SV[\Phi[{\bf v}],x]\right)=2g_{xx} \sin\theta(-\cos\theta\cos\phi +g_x\sin\theta).
$$
\item
$$
\frac{d}{dy}\left( SV[\Phi[{\bf v}],y]\right)=2g_{yy} \sin\theta(-\cos\theta\sin\phi +g_y\sin\theta).
$$
\end{enumerate}
\end{prop}
\demo
Direct calculations based on Definition \ref{def:sqvel} show the statements.
\enD

Furthermore, if we have two families of curves with respect to both $x$ and $y$ on a surface $M$, we also obtain the following theorems as byproducts of Propositions \ref{prop:curveturep}, \ref{prop:sv}.
\begin{thm}\label{thm:xy-curvature}
If $\Phi[{\bf v}]\colon U\to \R^2$ is  a one-parameter family of regular curves with respect to both $x$ and $y$,
then we have
\begin{eqnarray*}
\kappa[\Phi[{\bf v}],x]
\cdot
\kappa[\Phi[{\bf v}],y]
=-g_{xx}g_{yy}\frac{\sin\phi\cos\phi \sin^2\theta}{(\sin^2 \phi+P_{\theta,\phi}^2)^{\frac32}(\cos^2 \phi+Q_{\theta,\phi}^2)^{\frac32}}.
\end{eqnarray*}
\end{thm}

\begin{thm}\label{thm:xySV}
We have
\begin{eqnarray*}
\frac{d}{dx}\left( SV[\Phi[{\bf v}],x]\right)
\cdot
\frac{d}{dy}\left( SV[\Phi[{\bf v}],y]\right)
=4g_{xx}g_{yy}P_{\theta,\phi}Q_{\theta,\phi} \sin^2\theta.
\end{eqnarray*}
\end{thm}


\begin{rem}{\rm
Since the squared velocity can be defined for any smooth mappings,
we do not need any conditions to $\Phi[{\bf v}]$ in Theorem \ref{thm:xySV}.
}
\end{rem}

\section{Relationship between the signs of invariants}
In the d'Ocagne-Koenderink formula (\ref{KoenderinkFormula}),
the coincidence of the signs of the different quantities
is a remarkable feature,
which mathematically justifies the intuition that the approximate shape of an object can be determined by looking at its contour. 
Focusing on the sign of each quantity in the formulae in Theorems \ref{thm:y-zeta}, \ref{thm:x-zeta}, \ref{thm:xy-curvature} and \ref{thm:xySV},
we get the statements below.

\begin{cor}\label{cor:y-zeta}
Suppose that $\Phi[{\bf v}]\colon U\to \R^2$ is  a one-parameter family of regular curves with respect to $y$ and $sin \theta\not=0$.
Then we have
$$
\sign \kappa[\Phi[{\bf v}],y]
\cdot
\sign \frac{d}{dx}\left( SV[\Phi[{\bf v}],x]\right)
=
-\sign (g_{xx}g_{yy}) \cdot \sign P_{\theta,\phi}\cdot \sign \cos\phi.
$$
\end{cor}

\begin{cor}\label{cor:x-zeta}
Suppose that $\Phi[{\bf v}]\colon U\to \R^2$ is  a one-parameter family of regular curves with respect to $x$ and $sin \theta\not=0$.
Then we have
$$
\sign \kappa[\Phi[{\bf v}],x]
\cdot
\sign \frac{d}{dy}\left( SV[\Phi[{\bf v}],y]\right)
=
\sign (g_{xx}g_{yy}) \cdot \sign Q_{\theta,\phi}\cdot \sign \sin\phi.
$$
\end{cor}

\begin{cor}\label{cor:xy-curvature}
Suppose that $\Phi[{\bf v}]\colon U\to \R^2$ is  a one-parameter family of regular curves with respect to both $x$ and $y$ and $sin \theta\not=0$.
Then we have
\begin{eqnarray*}
\sign \kappa[\Phi[{\bf v}],x]
\cdot
\sign\kappa[\Phi[{\bf v}],y]
=-\sign (g_{xx}g_{yy}) \cdot \sign \sin2\phi. \end{eqnarray*}
\end{cor}
\begin{cor}\label{cor:xySV}
Suppose $sin \theta\not=0$.
Then we have
\begin{eqnarray*}
\sign \frac{d}{dx}\left( SV[\Phi[{\bf v}],x]\right)
\cdot
\sign \frac{d}{dy}\left( SV[\Phi[{\bf v}],y]\right)
=\sign (g_{xx}g_{yy})\cdot \sign P_{\theta,\phi}\cdot \sign Q_{\theta,\phi}.
\end{eqnarray*}
\end{cor}

Note that the following equation holds for a surface $M$ parametrized by $s(x,y)=(x,y,g(x,y))$:
$$
g_{xx}g_{yy}=K\left(1+g_x^2+g_y^2\right)^2+g_{xy}^2,
$$
where $K$ is the Gaussian curvature of $M$.
Thus, we have the following lemma.
\begin{lem}
\begin{enumerate}
    \item If $K(p)\ge0$ for $p\in M$, then $\sign K(p)=\sign( g_{xx}g_{yy}(p))$.
    \item 
    If $g_{xx} g_{yy}(p)\le0$ for $p\in M$, then $K(p)\le0$ holds for $p\in M$.
\end{enumerate}
\label{lem:gaussandg}
\end{lem}

\begin{rem}{\rm
     Combining Corollaries \ref{cor:y-zeta}--\ref{cor:xySV} and Lemma \ref{lem:gaussandg},
    it is possible to obtain relationships between the Gaussian curvature $K$ of a surface $M$ and invariants of the projections $\Phi[{\bf v}]$ as families of curves under appropriate assumptions. 
    }
\end{rem}


\section{Example}

Finally, as the simplest examples, we consider our formulae in a special case.
Suppose that a surface $M$ is locally parameterized by $s(x,y)=(x,y,g(x,y))^t$ 
with $g(x,y)$ having a critical point at the origin $0$,
that is, $g_x(0)=g_y(0)=0$.
In this case,
\begin{eqnarray*}
    P_{\theta,\varphi}(0)=-\cos\theta\cos\phi,\;\;
    Q_{\theta,\varphi}(0)=-\cos\theta\sin\phi
\end{eqnarray*}
hold.
Furthermore, we assume that the $x$-curve and $y$-curve of $M$ are along the principal directions at the origin, that is, $g_{xy}(0)=0$,
thus we have $g_{xx}g_{yy}(0)=K(0)$.
Then, from Corollaries \ref{cor:y-zeta}--\ref{cor:xySV},
we have simple relationships as in the following.

\begin{prop}\label{prop:y-zeta:simplecase}
Suppose $g_x(0)=g_y(0)=g_{xy}(0)=0$ and $\sin 2\theta\cos\phi\not=0$. Then $\Phi[{\bf v}]\colon U\to \R^2$ is a one-parameter family of regular curves with respect to $y$ around $0$,
and we have
$$
\sign \kappa[\Phi[{\bf v}],y](0)
\cdot
\sign \frac{d}{dx}\left( SV[\Phi[{\bf v}],x]\right)(0)
=
\sign K(0)\cdot \sign \cos\theta.
$$
\end{prop}

\begin{prop}\label{prop:x-zeta:simplecase}
Suppose $g_x(0)=g_y(0)=g_{xy}(0)=0$ and $\sin 2\theta\sin\phi\not=0$.
Then $\Phi[{\bf v}]\colon U\to \R^2$ is  a one-parameter family of regular curves with respect to $x$ around $0$,
and we have
$$
\sign \kappa[\Phi[{\bf v}],x](0)
\cdot
\sign \frac{d}{dy}\left( SV[\Phi[{\bf v}],y]\right)(0)
=
-\sign K(0)\cdot \sign \cos\theta.
$$
\end{prop}
\begin{prop}\label{prop:xy-curvature:simplecase}
Suppose $g_x(0)=g_y(0)=g_{xy}(0)=0$ and $\sin \theta\not=0$.
Then $\Phi[{\bf v}]\colon U\to \R^2$ is  a one-parameter family of regular curves with respect to both $x$ and $y$ around $0$,
and we have
\begin{eqnarray*}
\sign \kappa[\Phi[{\bf v}],x](0)
\cdot
\sign\kappa[\Phi[{\bf v}],y](0)
=-\sign K(0) \cdot \sign \sin2\phi. \end{eqnarray*}
\end{prop}
\begin{prop}\label{prop:xySV:simplecase}
Suppose $g_x(0)=g_y(0)=g_{xy}(0)=0$ and $\sin 2\theta\not=0$.
Then we have
\begin{eqnarray*}
\sign \frac{d}{dx}\left( SV[\Phi[{\bf v}],x]\right)(0)
\cdot
\sign \frac{d}{dy}\left( SV[\Phi[{\bf v}],y]\right)(0)
=\sign K(0)\cdot \sign \sin2\phi.
\end{eqnarray*}
\end{prop}

\begin{exam}\label{ex:ellip}
    Let $M$ be a surface expressed as the graph of a function $g(x,y)=-x^2- y^2$ on $\R^2$,
    where $K(x,y)>0$ holds for any $(x,y)\in\mathbb{R}^2$. 
    Take the view direction 
    ${\bf v}=(\sin \theta \cos \phi,\sin \theta \sin \phi, \cos \theta)^t\in S^2$ so that $\sin 2\theta\cos\phi\not=0$.
    Then $\Phi[{\bf v}](x,y)$ is a one-parameter family of regular curves with respect to $y$ around $0$.
    In this case,
    $$
    \kappa[\Phi[{\bf v}],y](0)=\frac{2\cos\phi\sin\theta}{(\cos^2\phi+\cos^2\theta\sin^2\phi)^{\frac32}},\;\; \frac{d}{dx}\left( SV[\Phi[{\bf v}],x]\right)(0)=2\sin2\theta\cos\phi.
    $$
    Figure \ref{fig:ellipproj} represents the case with $\phi=0$ and $\frac\pi2<\theta<\pi$,
    where $\kappa[\Phi[{\bf v}],y](0)>0$ and $\frac{d}{dx}\left( SV[\Phi[{\bf v}],x]\right)(0)<0$.
    Thus we see that the relationship in Proposition \ref{prop:y-zeta:simplecase} holds in this case.

    \end{exam}

\begin{figure}
\centering
\includegraphics[clip,width=12.0cm]{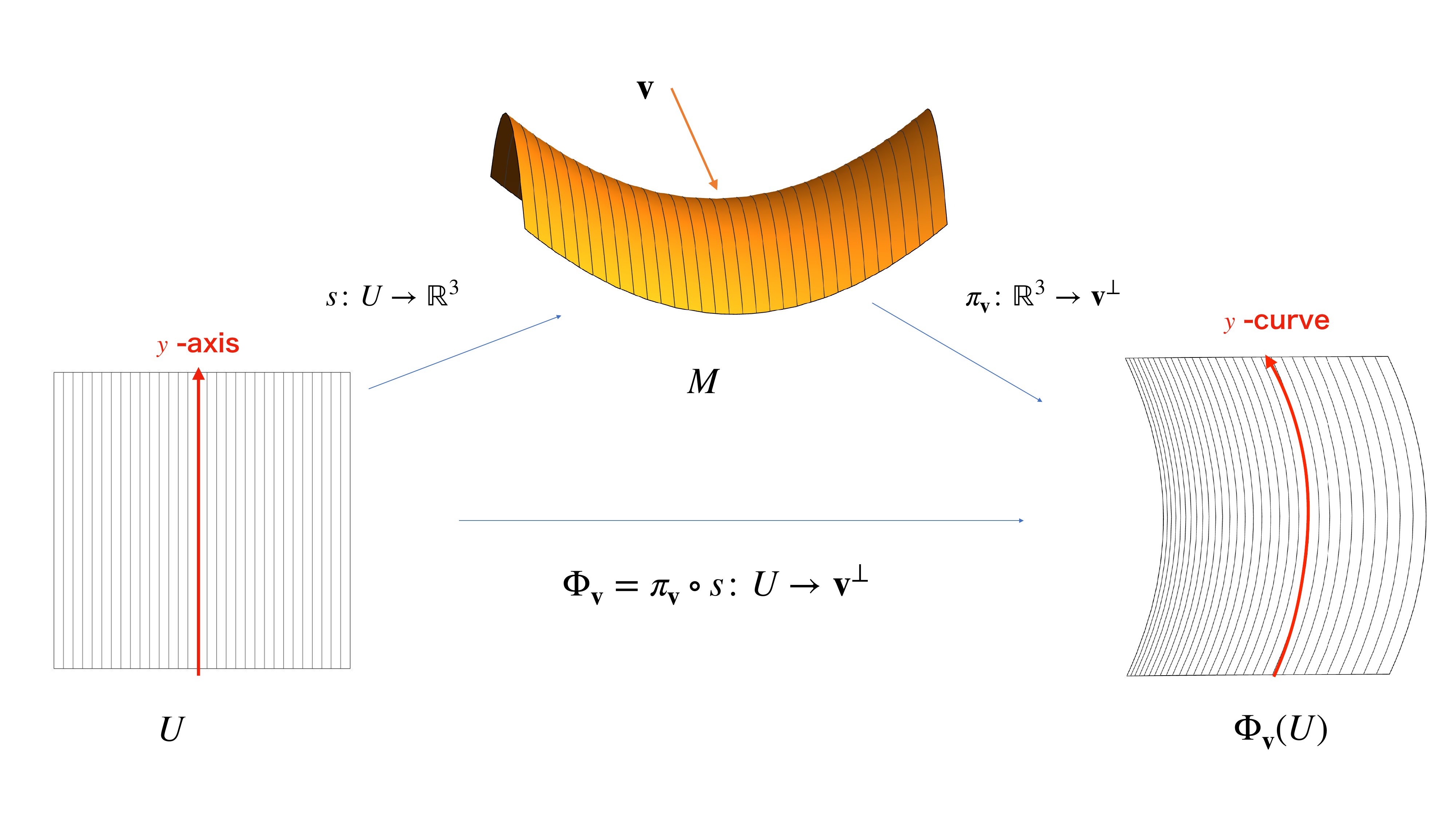}
\caption{
$M$ is hyperbolic, that is, $K>0$. Here $\kappa[\Phi[{\bf v}],y]>0$ and
$\frac{d}{dx}\left( SV[\Phi[{\bf v}],x]\right)<0$ hold for the one-parameter family of regular curves $\Phi[{\bf v}]$ with respect to $y$.
}
\label{fig:hypproj}
\end{figure}

\begin{exam}\label{ex:hyp}
    Let $M$ be a surface expressed as the graph of a function $g(x,y)=x^2- y^2$ on $\R^2$,
    where $K(x,y)<0$ holds for any $(x,y)\in\mathbb{R}^2$. 
    Take the view direction 
    ${\bf v}=(\sin \theta \cos \phi,\sin \theta \sin \phi, \cos \theta)^t\in S^2$ so that $\sin 2\theta\cos\phi\not=0$.  
    Then $\Phi[{\bf v}](x,y)$ is a one-parameter family of regular curves with respect to $y$ around $0$.
    In this case, it is easy to see that
    $$
    \kappa[\Phi[{\bf v}],y](0)=\frac{2\cos\phi\sin\theta}{(\cos^2\phi+\cos^2\theta\sin^2\phi)^{\frac32}},\;\; \frac{d}{dx}\left( SV[\Phi[{\bf v}],x]\right)(0)=-2\sin2\theta\cos\phi.
    $$
    Figure \ref{fig:hypproj} represents the case with $\phi=0$ and $\frac\pi2<\theta<\pi$,
    where $\kappa[\Phi[{\bf v}],y](0)>0$ and $\frac{d}{dx}\left( SV[\Phi[{\bf v}],x]\right)(0)>0$.
        Thus we see that the relationship in Proposition \ref{prop:x-zeta:simplecase} holds in this case.

    \end{exam}

    \begin{figure}
\centering
\includegraphics[clip,width=12.0cm]{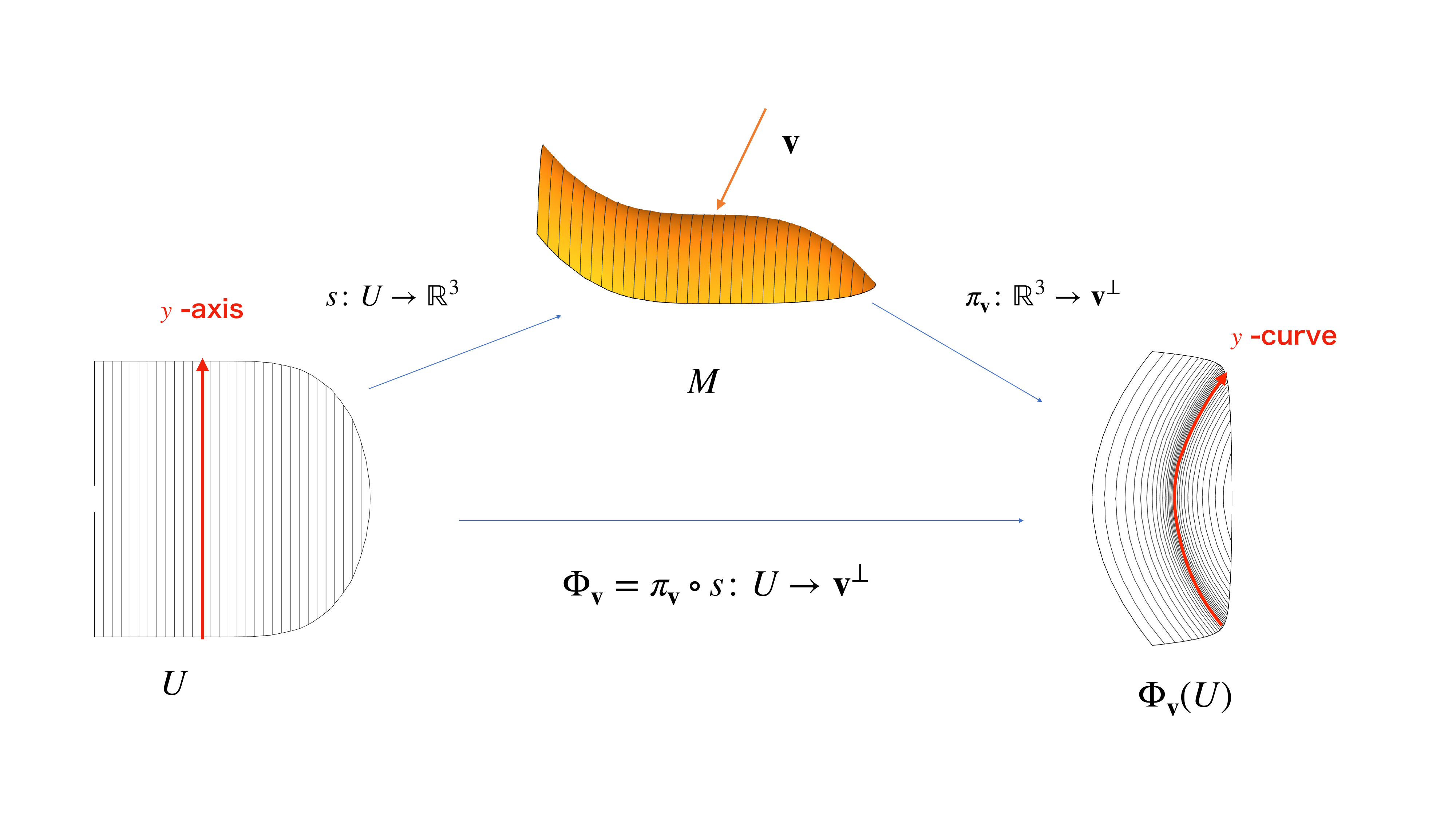}
\caption{
$M$ contains parabolic points $p$ where $K(p)=0$. $\frac{d}{dx}\left( SV[\Phi[{\bf v}],x]\right)(p)=0$ holds for the one-parameter family of regular curves $\Phi[{\bf v}]$ with respect to $y$.
}
\label{fig:paraboproj}
\end{figure}

    \begin{exam}\label{ex:para}
    Let $M$ be a surface expressed as the graph of a function $g(x,y)=y^2-x^3$ on $\R^2$, where $K(0)=0$ holds. 
    Take the view direction 
    ${\bf v}=(\sin \theta \cos \phi,\sin \theta \sin \phi, \cos \theta)^t\in S^2$ so that $\sin 2\theta\cos\phi\not=0$.
    Then $\Phi[{\bf v}](x,y)$ is a one-parameter family of regular curves with respect to $y$ around $0$.
    In this case, it is easy to see that
    $$
    \frac{d}{dx}\left( SV[\Phi[{\bf v}],x]\right)(0)=0.
    $$
    Thus we see that the relationship in Proposition \ref{prop:y-zeta:simplecase} holds in this case.
    Figure \ref{fig:paraboproj} represents the case.
\end{exam}

    Observe that, in Examples \ref{ex:ellip} -- \ref{ex:para}, the difference in the patterns of the projection images corresponds to the difference in the signs of the Gaussian curvatures of the surfaces.
    Compare also Figures \ref{fig:ellipproj} -- \ref{fig:paraboproj}.

\begin{exam}
    Let $M$ be a surface expressed as the graph of a function $g(x,y)=-x^2- y^2$ on $\R^2$,
    where $K(x,y)>0$ holds for any $(x,y)\in\mathbb{R}^2$. 
    Take the view direction 
    ${\bf v}=(\sin \theta \cos \phi,\sin \theta \sin \phi, \cos \theta)^t\in S^2$ so that $\sin 2\theta\not=0$. Then $\Phi[{\bf v}]$ is a one-parameter family of regular curves with respect to both $x$ and $y$ around $0$. We can see
    $$
    \kappa[\Phi[{\bf v}],x](0)=-\frac{2\sin\phi\sin\theta}{\sin^2\phi+\cos^2\theta\cos^2\phi},\; 
    \kappa[\Phi[{\bf v}],y](0)=-\frac{2\cos\phi\sin\theta}{\cos^2\phi+\cos^2\theta\sin^2\phi},\;\;     
    $$
    and
    $$
    \frac{d}{dx}\left( SV[\Phi[{\bf v}],x]\right)(0)=2\sin2\theta\cos\phi,\;
\frac{d}{dy}\left( SV[\Phi[{\bf v}],y]\right)(0)=2\sin2\theta\sin\phi.
    $$
    See also Figures \ref{fig:ellipprojcurvaturexy}, \ref{fig:ellipprojSVxy} where $\frac\pi2<\theta<\pi$ and $\pi<\phi<\frac32 \pi$.
    Here the relationships in Propositions \ref{prop:xy-curvature:simplecase} and \ref{prop:xySV:simplecase} hold.

\begin{figure}
\centering
\includegraphics[clip, width=13.0cm]{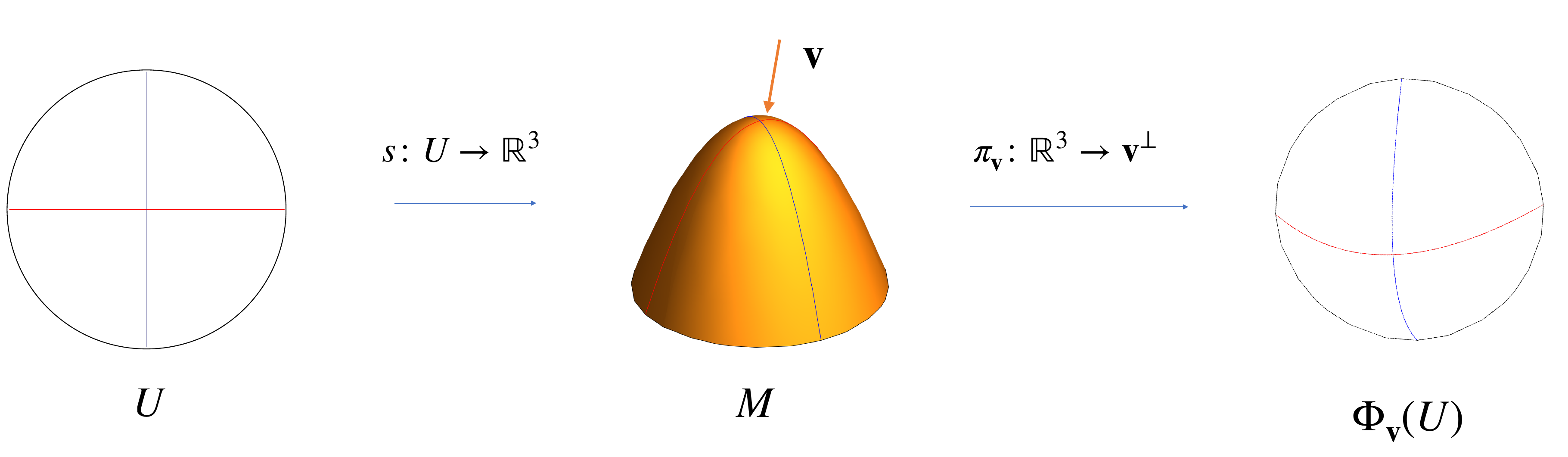}
\caption{
The surface $M$ is elliptic, that is, $K>0$, and $\kappa[\Phi[{\bf v}],y](0)<0$ and $\kappa[\Phi[{\bf v}],x](0)>0$
hold for the one-parameter family of regular curves $\Phi[{\bf v}]$ with respect to $x$ and $y$.
}
\label{fig:ellipprojcurvaturexy}
\end{figure}

\begin{figure}
\centering
\includegraphics[clip, width=13.0cm]{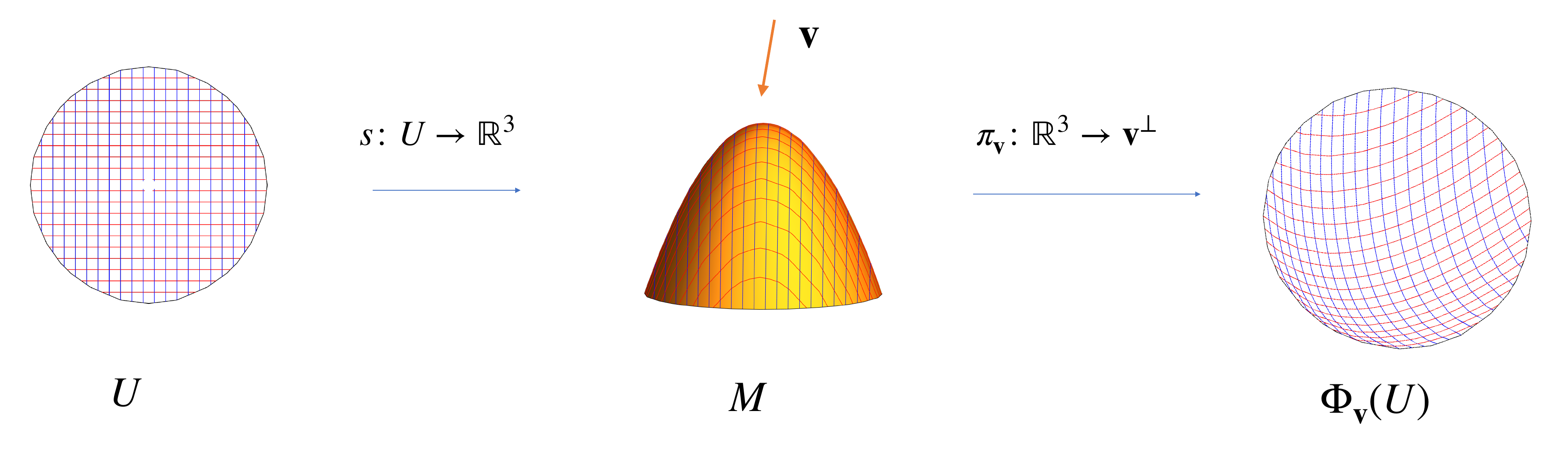}
\caption{
The surface $M$ is elliptic, that is, $K>0$, and 
$
\frac{d}{dx}\left( SV[\Phi[{\bf v}],x]\right)(0)
>0$ $
\frac{d}{dy}\left( SV[\Phi[{\bf v}],y]\right)(0)
>0$
hold for the one-parameter family of regular curves $\Phi[{\bf v}]$ with respect to $x$ and $y$.
}
\label{fig:ellipprojSVxy}
\end{figure}

\end{exam}

\noindent
Ken Anjyo, 
\\
OLM Digital Inc., Japan and IMAGICA GROUP Inc., Japan\\
\\
\\
Yutaro Kabata,
\\
School of Information and Data Sciences, 
Nagasaki University, 
Nagasaki 852-8131, Japan.\\
E-mail address: kabata@nagasaki-u.ac.jp
\\

\end{document}